\documentclass[11pt]{article}
\usepackage{cite}
\usepackage{mathrsfs}
\usepackage{amsfonts}
\usepackage{amsmath}
\usepackage{amsfonts,amssymb,color}
\usepackage{dsfont}
\usepackage{curves}
\usepackage{mathrsfs}
\usepackage{pifont}
\usepackage{amssymb}
\usepackage{latexsym,amsmath,amssymb,amsfonts,epsfig,graphicx,cite,psfrag}
\usepackage{eepic,color,colordvi,amscd}

\newtheorem{theorem}{Theorem}[section]

\newtheorem{lemma}{Lemma}[section]

\newtheorem{claim}{Claim}[section]
\newtheorem{conjecture}{Conjecture}[section]

\newcommand{\qed}{\hfill\rule{0.5em}{0.809em}}

\def\emptyset{\mbox{{\rm \O}}}
\def\bar{\overline}

\def\qed{\hfill \rule{4pt}{7pt}}

\def\pf{\noindent {\it Proof. }}

\begin{document}
\title{A tight linear chromatic bound for ($P_3\cup P_2, W_4$)-free graphs}

\author{{Rui Li $^{a,}$\footnote{Email address: lirui@hhu.edu.cn}\quad Jinfeng Li $^{a,}$\footnote{Email address: 1345770246@qq.com}\quad  Di Wu$^{b,}$\footnote{Email address: 1975335772@qq.com
}}\\
       {\small $^a$ School of Mathematics, Hohai University}\\
        {\small 8 West Focheng Road, Nanjing, 211100, China}\\
        {\small $^b$ School of Mathematical Science, Nanjing Normal University}\\
        {\small 1 Wenyuan Road, Nanjing, 210046, China}
        }
\date{}
\maketitle

\begin{abstract}
 
For two vertex disjoint graphs $H$ and $F$,  we use  $H\cup F$ to denote  the graph with vertex set $V(H)\cup V(F)$ and edge set $E(H)\cup E(F)$,  and use $H+F$ to denote the graph with vertex set $V(H)\cup V(F)$ and edge set $E(H)\cup E(F)\cup\{xy\;|\; x\in V(H), y\in V(F)$$\}$. A $W_4$ is the graph $K_1+C_4$. In this paper, we prove that  $\chi(G)\le 2\omega(G)$ if $G$ is a ($P_3\cup P_2, W_4$)-free graph. This bound is tight when $\omega =2$ and $3$, and improves the main result  of Wang and Zhang \cite{WZ22}. Also, this bound partially generalizes some results of Prashant {\em et al.}\cite{PA22}.

\begin{flushleft} {\em Key words and
phrases:}  chromatic number; clique number; $\chi$-binding function; $(P_3\cup P_2)$-free graphs 

\end{flushleft}
\end{abstract}

\section{Introduction}

All graphs considered in this paper are finite and simple. We use $P_k$ and $C_k$ to denote a {\em path} and a {\em cycle} on $k$ vertices respectively, and follow \cite{BM1976} for undefined notations and terminology. Let $G$ be a graph, and $X$ be a subset of $V(G)$. We use $G[X]$ to denote the subgraph of $G$ induced by $X$, and call $X$ a {\em clique} ({\em independent set}) if $G[X]$ is a complete graph (has no edge). The {\em clique number} $\omega(G)$ of $G$ is the maximum size taken over all cliques of $G$.

For $v\in V(G)$, let $N_G(v)$ be the set of vertices adjacent to $v$, $d_G(v)=|N_G(v)|$, $N_G[v]=N_G(v)\cup \{v\}$, $M_G(v)=V(G)\setminus N_G[v]$, $M_G[v]=V(G)\setminus N_G(v)$. For $X\subseteq V(G)$, let $N_G(X)=\{u\in V(G)\setminus X\;|\; u$ has a neighbor in $X\}$ and $M_G(X)=V(G)\setminus (X\cup N_G(X))$, $N_G[X]=N_G(X)\cup X$ and $M_G[X]=M_G(X)\cup X$. If it does not cause any confusion, we will omit the subscript $G$ and simply write $N(v), d(v), N[v], M(v), M[v], N(X), N[X], M(X)$ and $M[X]$.  Let $\Delta(G)$ ($\delta(G)$) denote the maximum (minimum) degree of $G$.

Let $G$ and $H$ be two vertex disjoint graphs. The {\em union} $G\cup H$ is the graph with $V(G\cup H)=V(G)\cup (H)$ and $E(G\cup H)=E(G)\cup E(H)$. The union of $k$ copies of the same graph $G$ will be denoted by $kG$. The {\em join} $G+H$ is the graph with $V(G+H)=V(G)\cup V(H)$ and $E(G+H)=E(G)\cup E(H)\cup\{xy\;|\; x\in V(G), y\in V(H)$$\}$. The complement of a graph $G$ will be denoted by $\bar G$.  We say that $G$ induces $H$ if $G$ has an induced subgraph isomorphic to $H$, and say that $G$ is $H$-free otherwise. Analogously, for a family $\cal H$ of graphs, we say that $G$ is ${\cal H}$-free if $G$ induces no member of ${\cal H}$.

Let $k$ be a positive integer, and let $[k]=\{1, 2, \ldots, k\}$. A $k$-{\em coloring} of $G$ is a mapping $c: V(G)\mapsto [k]$ such that $c(u)\neq c(v)$ whenever $u\sim v$ in $G$. The {\em chromatic number} $\chi(G)$ of $G$ is the minimum integer $k$ such that $G$ admits a $k$-coloring.  It is certain that $\chi(G)\ge \omega(G)$. A {\em perfect graph} is one such that  $\chi(H)=\omega(H)$ for all of its induced subgraphs $H$.  A family $\cal G$ of graphs is said to be $\chi$-{\em bounded} if there is a function $f$ such that $\chi(G)\le f(\omega(G))$ for every $G\in\cal G$, and if such a function does exist for $\cal G$, then $f$ is said to be a {\em binding function} of $\cal G$ \cite{Gy75}.

Erd\"{o}s \cite{ER59} proved that for any positive integers $k,l\ge3$, there exists a graph $G$ with $\chi(G)\ge k$ and no cycles of length less than $l$. This result motivates us to study the chromatic number of $F$-free graphs, where $F$ is a forest (a disjoint union of trees). Gy\'{a}rf\'{a}s \cite{Gy75} and Sumner \cite{Su81} independently, proposed the following famous conjecture.

\begin{conjecture}\label{tree}\cite{Gy75, Su81}
	Let $F$ be a forest. Then $F$-free graphs are $\chi$-bounded.
\end{conjecture}

The classes of $2K_2$-free graphs and $P_5$-free graphs have attracted a great deal of interest in recent years.  Up to now, the best known $\chi$-binding function for $2K_2$-free graphs is $f(\omega)=\binom{\omega+1}{2}-2\lfloor \frac{\omega}{3} \rfloor$ \cite{GM22}, and the best known $\chi$-binding function for $P_5$-free graphs is $f(\omega)=\omega(G)^{\log_2\omega(G)}$ \cite{SS21}. We refer the interested readers to \cite{GM22} for results of $P_5$-free graphs, and to \cite{BBS19,KM18,PA23,Wa80}  for results of $2K_2$-free graphs, and to \cite{RS04,SR19,SS20} for more results and problems about the $\chi$-bounded problem. Since $(P_3\cup P_2)$ is a forest with five vertices and the class of $(P_3\cup P_2)$-free graphs is a superclass of $2K_2$-free graphs, many scholars began to show interest in $(P_3\cup P_2)$-free graphs. The best known $\chi$-binding function for $(P_3\cup P_2)$-free graphs is $f(\omega)=\frac{1}{6}\omega(\omega+1)(\omega+2)$ \cite{BA18}.  Fortunately, this function is polynomial, but we would prefer a linear one, even for some subclasses of $(P_3\cup P_2)$-free graphs.

\begin{figure}[htbp]\label{fig-1}
	\begin{center}
		\includegraphics[width=12cm]{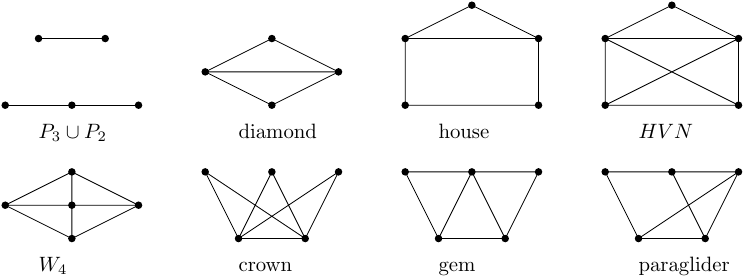}
	\end{center}
	\vskip -15pt
	\caption{Illustration of $P_3\cup P_2$ and some forbidden configurations.}
\end{figure}

A {\em diamond} is the graph $K_1+P_3$, a {\em house} is just the complement of $P_5$, an $HVN$ is a $K_4$ together with one more vertex which is adjacent to exactly two vertices of $K_4$, a $W_i$ is the graph $K_1+C_i$ ($i$ is an integer and $i\ge4$),  a {\em crown} is the graph $K_1+K_{1,3}$, a {\em gem}  is the graph $K_1+P_4$, and  a {\em paraglider} is the graph obtained from a diamond by adding a vertex joining to its two vertices of degree 2 (See Figure \ref{fig-1}).

In \cite{PA22}, Prashant {\em  et al.} proved that if $G$ is ($P_3\cup P_2$, diamond)-free, then $\chi(G)\le4$ when $\omega(G)=2$, $\chi(G)\le6$ when $\omega(G)=3$, $\chi(G)\le4$ when $\omega(G)=4$, and $G$ is perfect when $\omega(G)\ge5$, and they also proved \cite{PA22} that $\chi(G)\le\omega(G)+1$ if $G$ is a $(P_3\cup P_2, HVN)$-free graph with $\omega(G)\ge4$.  As a superclass of ($P_3\cup P_2$, diamond)-free graphs, Cameron {\em  et al.} \cite{CHM21} proved that $\chi(G)\le\omega(G)+3$ if $G$ is ($P_6$, diamond)-free, this bound is optimal.  Very Recently, Wu and Xu \cite{WX23} proved that $\chi(G)\le\frac{1}{2}\omega^2(G)+\frac{3}{2}\omega(G)+1$ if $G$ is ($P_3\cup P_2$, crown)-free, Char and Karthick \cite{CK22} proved that $\chi(G)\le$ max $\{\omega(G)+3, \lfloor \frac{3\omega(G)}{2} \rfloor-1\}$ if $G$ is a ($P_3\cup P_2$, paraglider)-free graph with $\omega(G)\ge3$, Prashant {\em et al.} proved that $\chi(G)\le2\omega(G)$ if $G$ is ($P_3\cup P_2$, gem)-free, and  Li {\em et al.}\cite{LLW23} proved that $\chi(G)\le2\omega(G)$ if $G$ is ($P_3\cup P_2$, house)-free.

In \cite{WZ22}, Wang and Zhang proved that if $G$ is a $(P_3\cup P_2, K_3)$-free
graph, then $\chi(G)\le3$ unless $G$ is one of eight graphs with $\Delta(G)=5$ and $\chi(G)=4$, and they proved \cite{WZ22} $\chi(G)\le3\omega(G)$ if $G$ is $(P_3\cup P_2, W_4)$-free. 

In this paper, we prove that

\begin{theorem}\label{W4}
	$\chi(G)\le2\omega(G)$ if $G$ is $(P_3\cup P_2, W_4)$-free.
\end{theorem}

This generalizes some results of Wang and Zhang \cite{WZ22}. In particular, $\chi(G)\le6$ if $G$ is $(P_3\cup P_2, W_4)$-free with $\omega(G)=3$, which generalizes the result of Prashant {\em et al.}\cite{PA22} that $\chi(G)\le6$ if $G$ is ($P_3\cup P_2$, diamond)-free with $\omega(G)=3$.

\medskip

Let $H$ be the Mycielski-Gr\"{o}stzsch graph (see Figure 2), and $H'$ be the complement of Schl\"{a}fli graph (see https://houseofgraphs.org/graphs/19273). Notice that $H$ and $H'$ are both $(P_3\cup P_2, W_4)$-free, and $\omega(H)=2,\chi(H)=4$, and $\omega(H')=3,\chi(H')=6$. Therefore, the $\chi$-binding function $f=2\omega$ is tight when $\omega=2$ and $3$.

\begin{figure}[htbp]\label{fig-2}
	\begin{center}
		\includegraphics[width=4cm]{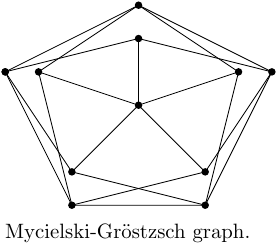}
	\end{center}
	\vskip -15pt
	\caption{Mycielski-Gr\"{o}stzsch graph.}
\end{figure}

\section{Preliminary and notations}

Let $G$ be a graph, $v\in V(G)$, and let $X$ and $Y$ be two subsets of $V(G)$. We say that $v$ is {\em complete} to $X$ if $v$ is adjacent to all vertices of $X$, and say that $v$ is {\em anticomplete} to $X$ if $v$ is not adjacent to any vertex of $X$. We say that $X$ is complete (resp. anticomplete) to $Y$ if each vertex of $X$ is complete (resp. anticomplete) to $Y$. For $u, v\in V(G)$, we simply write $u\sim v$ if $uv\in E(G)$, and write $u\not\sim v$ if $uv\not\in E(G)$.

A {\em hole} of $G$ is an induced cycle of length at least 4, and a {\em $k$-hole} is a hole of length $k$. A $k$-hole is called an {\em odd hole} if $k$ is odd, and is called an {\em even hole} otherwise. An {\em antihole} is the complement of some hole. An odd (resp. even) antihole is defined analogously. The famous {\em Strong Perfect Graph Theorem} states that

\begin{theorem}\label{Perfect}\cite{CRSR06}
	A graph is perfect if and only if it induces neither an odd hole nor an odd antihole.
\end{theorem}

The following lemmas  and theorems will be used in our proof.

\begin{lemma}\label{C3}\cite{WZ22}
	If $G$ is a $(P_3\cup P_2, C_3)$-free  graph, then $\chi(G)\le 4$.
\end{lemma}

\begin{lemma}\label{C4}\cite{CK10}
	If $G$ is a $(P_3\cup P_2, C_4)$-free  graph, then $\chi(G)\le \lceil \frac{5\omega(G)}{4} \rceil$.
\end{lemma}

\begin{theorem}\label{gem}\cite{PAM23}
	If $G$ is a $($$P_3\cup P_2$, gem$)$-free graph, then $\chi(G)\le2\omega(G)$.
\end{theorem}

\begin{theorem}\label{K3}(1.8 of \cite{ CS10})
	Let $G$ be a $($$K_1\cup K_3$$)$-free graph. If $G$ contains a $K_3$, then $\chi(G)\le2\omega(G)$.
\end{theorem}

Let $G$ be a ($P_3\cup P_2, W_4$)-free graph. By Lemma \ref{C3}, we have that $\chi(G)\le 4=2\omega(G)$ if $\omega(G)=2$. Therefore, we may divide the proof process of Theorem \ref{W4} into three sections depending on $\omega(G)\ge 5$ or $G$ is $K_5$-free, or $G$ is $K_4$-free. 

\section{$\omega\ge5$}

In this section, we consider ($P_3\cup P_2, W_4$)-free graphs with $\omega\ge5$. 

\begin{lemma}\label{ge5}
	If $G$ is a $(P_3\cup P_2, W_4)$-free graph with $\omega(G)\ge5$, then $\chi(G)\le2\omega(G)$.
\end{lemma}

\pf Let $G$ be a ($P_3\cup P_2, W_4$)-free graph with $\omega(G)\ge5$, we will complete the proof of Lemma \ref{ge5} by the following two claims. 

\begin{claim}\label{c1}
	If $G$ contains an induced $W_5$, then $\chi(G)\le2\omega(G)$.
\end{claim}

\pf Let $Q$ be an induced $W_5$ in $G$ with $V(Q)=\{u,v_1,v_2,v_3,v_4,v_5\}$ such that $d_Q(u)=5$ and $C=v_1v_2v_3v_4v_5v_1$ is an induced $C_5$. Since $G$ is $W_4$-free, we have that $G[N(u)]$ is $C_4$-free, which implies that $\chi(G[N(u)])\le \lceil \frac{5}{4}(\omega(G)-1) \rceil$ by Lemma \ref{C4}. Now, we only need to color $G[M(u)]$.

Let $x\in M(u)$. Suppose $G[N_C(x)]$ contains an induced $P_3$, say $v_1v_2v_3$ by symmetry. Then $\{x,v_1,v_2,v_3,u\}$ induces a $W_4$, a contradiction. So, $G[N_C(x)]$ is $P_3$-free. Let $M_i
=\{x\in M(u)|x$ is anticomplete to $\{v_i,v_{i+2}\}\}$ and $M=\cup^5_{i=1}M_i$ (The summation of subindex is taken
modulo $5$, and we set $5+1\equiv 1$). Since $G[N_C(x)]$ is $P_3$-free, it follows that $x\in M$. Moreover, $M_i$ is independent for each $i\in\{1,2,3,4,5\}$. Otherwise, suppose $M_1$ is not independent that it has two adjacent vertices $x_1$ and $x_2$. Then $\{v_1,u,v_3,x_1,x_2\}$ induces a $P_3\cup P_2$, a contradiction. So,  $\chi(G[M(u)])\le5$.

Therefore, $\chi(G)\le\chi(G[N(u)])+\chi(G[M[u]])\le \lceil \frac{5}{4}(\omega(G)-1) \rceil+5\le2\omega(G)$ as $\omega(G)\ge5$. This proves Claim \ref{c1}.  \qed

\begin{claim}\label{c2}
	If $G$ is $W_5$-free and contains an induced $P_2\cup K_3$, then $\chi(G)\le2\omega(G)$.
\end{claim}

\pf Let $Q$ be an induced $P_2\cup K_3$ in $G$ with $\{u_1,u_2\}\subseteq V(Q)$ such that $d_Q(u_1)=d_Q(u_2)=1$. Let $N=N(u_2)\setminus N[u_1]$ and $M=M(\{u_1,u_2\})$. Since $G$ is $(P_3\cup P_2, W_4, W_5)$-free, it follows that $G[N(u_i)]$ is $(P_3\cup P_2, C_4, C_5)$-free, for each $i\in\{1,2\}$. By the famous Strong Perfect Theorem, we have that $G[N(u_1)]$ and $G[N(u_2)]$ are both perfect. Particularly, $G[N]$ is perfect. Moreover, $G[M]$ is $P_3$-free as $G$ is $(P_3\cup P_2)$-free, which implies that $G[M]$ is perfect.

Choose a maximum clique $C_1$ in $N$, and choose a maximum clique $C_2$ in $M$. Note that $|C_2|\ge3$ as the $Q$ exists. For each vertex $u\in C_1$, if $|N(u)\cap C_2|\le |C_2|-2$, let $v,v'\in C_2$ such that $u\not\sim v, u\not\sim v'$, then $\{u_1,u_2,u,v,v'\}$ induces a $P_3\cup P_2$, a contradiction. So, we have that

\begin{equation}\label{eqa-1}
	\mbox{for each $u\in C_1$, $|N(u)\cap C_2|\ge |C_2|-1$}.
\end{equation}

Suppose there exist two vertices $x$ and $y$ in $C_1$ such that $N(x)\cap C_2\not\subseteq N(y)\cap C_2$ and $N(y)\cap C_2\not\subseteq N(x)\cap C_2$. Then there exist two vertices $v_1$ and $v_2$ such that $x\sim v_1, x\not\sim v_2, y\sim v_2, y\not\sim v_1$. Moreover, by (\ref{eqa-1}) and $|C_2|\ge3$, there exists a vertex $v_3\in C_2$ such that $v_3\sim x$ and $v_3\sim y$. Now, $\{v_3,x,v_1,v_2,y\}$ induces a $W_4$, a contradiction. Therefore, for any two vertices $x$ and $y$ in $C_1$, either $N(x)\cap C_2\subseteq N(y)\cap C_2$ or $N(y)\cap C_2\subseteq N(x)\cap C_2$, which implies that $|C_1|+|C_2|\le\omega(G)+1$ by (\ref{eqa-1}). And thus, $\chi(G[N\cup M\cup \{u_1,u_2\}])\le\chi(G[N])+\chi(G[M\cup\{u_1,u_2\}])\le\omega(G[N])+\omega(G[M])\le\omega(G)+1$.

So, $\chi(G)\le \chi(G[N(u_1)])+\chi(G[N\cup M\cup \{u_1,u_2\}])\le\omega(G)-1+\omega(G)+1=2\omega(G)$. This proves Claim \ref{c2}.\qed

\medskip

By Claim \ref{c1} and \ref{c2}, we may assume that $G$ is $(P_3\cup P_2, W_4, W_5, P_2\cup K_3)$-free. Let $u_1u_2$ be an edge of $G$. Then $G[M(u_1,u_2)]$ is $(P_3, K_3)$-free, which implies that $\chi(G[M(u_1,u_2)])\le2$. Furthermore, $G[N(u_1)]$ and $G[N(u_2)]$ are both perfect. So, $\chi(G)\le \chi(G[N(u_1)]) +\chi(G[N(u_2)])+\chi(G[M(u_1,u_2)])\le\omega(G)-1+\omega(G)-1+2=2\omega(G)$.

This proves Lemma \ref{ge5}.\qed

\section{$K_5$-free}

In this section, we consider ($P_3\cup P_2, W_4, K_5$)-free graphs. 

\begin{lemma}\label{=4}
	If $G$ is a $(P_3\cup P_2, W_4, K_5)$-free graph, then $\chi(G)\le8$.
\end{lemma}
\pf Let $G$ be a $(P_3\cup P_2, W_4, K_5)$-free graph, we will complete the proof of Lemma \ref{=4} by the following five claims. 

\begin{claim}\label{c3}
	If $G$ contains an induced $2K_3$, then $\chi(G)\le8$.
\end{claim}
\pf Let $G[\{v_1,v_2,v_3,u_1,u_2,u_3\}]$ be an induced $2K_3$ in $G$ such that $v_1v_2v_3v_1$ is a triangle.
 Since $G$ is ($P_3\cup P_2$)-free, we have that $M(\{v_1,v_2,v_3\})$ is the union of cliques. Without loss of generality, we suppose that $\{u_1,u_2,u_3\}\subseteq C$, where $C$ is the maximum clique of $G[M(\{v_1,v_2,v_3\})]$. Let $N=\{v\in N(\{v_1,v_2,v_3\})| |N(v)\cap \{v_1,v_2,v_3\}|\le2\}$ and $x\in N$. 

Without loss of generality, $x\sim v_1$ and $x\not\sim v_2$. If $|\{u_1,u_2,u_3\}\cap N(x)|\le1$, we may by symmetry assume that $x\not\sim u_1$ and $x\not\sim u_2$, then $\{x,v_1,v_2,u_1,u_2\}$ induces a $P_3\cup P_2$, a contradiction. So, $|N(x)\cap \{u_1,u_2,u_3\}|\ge2$. Let $N_1=\{v\in N | |N(v)\cap\{u_1,u_2,u_3\}|=2\}$ and $N_2=\{v\in N | |N(v)\cap\{u_1,u_2,u_3\}|=3\}$. It is obvious that $N=N_1\cup N_2$ and $N_2$ is independent as $G$ is $K_5$-free. 

Let $y\in N_1$. Then we may assume that $N(y)\cap\{u_1,u_2,u_3\}=\{u_1,u_2\}$. Suppose there exists a vertex $y'\in N_1$ such that $y\sim y'$. If $N(y')\cap\{u_1,u_2,u_3\}\ne\{u_1,u_2\}$, let $N(y')\cap\{u_1,u_2,u_3\}=\{u_1,u_3\}$, then $\{y,y',u_1,u_2,u_3\}$ induces a $W_4$, a contradiction. So, for any adjacent vertices $y$ and $y'$ in $N_1$, $N(y)\cap\{u_1,u_2,u_3\}=N(y')\cap\{u_1,u_2,u_3\}$. Therefore, $G[N_1]$ is $(K_3, C_4)$-free as $G$ is $(K_5, W_4)$-free. By Lemma \ref{c2}, we have that $\chi(G[N_1])\le3$.

If $|C|=3$, that is to say, $C=\{u_1,u_2,u_3\}$, then we can directly deduce that $\chi(G[C\cup N_2\cup\{v_1,v_2,v_3\}])\le4$. Otherwise $|C|=4$, let $C=\{u_1,u_2,u_3,u'\}$, then $u'$ is anticomplete to $N_2$ as $G$ is $K_5$-free, which still implies that $\chi(G[C\cup N_2\cup\{v_1,v_2,v_3\}])\le4$.

So, $\chi(G[M[\{v_1,v_2,v_3\}]\cup N_2])\le4$.  Notice that $N(\{v_1,v_2,v_3\})=(N(v_1)\cap N(v_2)\cap N(v_3))\cup N$ and $N=N_1\cup N_2$. Since $G$ is $K_5$-free, it follows that $N(v_1)\cap N(v_2)\cap N(v_3)$ is independent. So, $\chi(G[N(\{v_1,v_2,v_3\})\setminus N_2])\le4$. Now, $\chi(G)\le \chi(G[M[\{v_1,v_2,v_3\}]])+\chi(G[N(\{v_1,v_2,v_3\})])\le\chi(G[M[\{v_1,v_2,v_3\}]\cup N_2]) +\chi(G[N(\{v_1,v_2,v_3\})\setminus N_2])\le4+4=8$. This proves Claim \ref{c3}. \qed

\begin{claim}\label{c4}
	If $G$ is $2K_3$-free and contains an induced $P_2\cup K_4$, then $\chi(G)\le8$.
\end{claim}
\pf Let $G[\{v_1,v_2,u_1,u_2,u_3,u_4\}]$ be an induced $P_2\cup K_4$ in $G$ such that $u_1u_2u_3u_4u_1$ is a $K_4$. Let $N_1=N(v_1)\setminus N[v_2], N_2=N(v_2)\setminus N[v_1], N_3=N(v_1)\cap N(v_2)$, and $N_4=\{v\in N_3| |N(v)\cap \{u_1,u_2,u_3,u_4\}|=3\}$. Note that $N(\{v_1,v_2\})=N_1\cup N_2\cup N_3$. We will prove that 

\begin{equation}\label{eqa-3}
	\mbox{$N_1\cup N_2\cup N_4$ is independent}.
\end{equation}

Let $x$ be any vertex in $N_1$. If $|\{u_1,u_2,u_3,u_4\}\setminus N(x)|\ge2$, let $x\not\sim u_1$ and $x\not\sim u_2$, then $\{x,v_1,v_2,u_1,u_2\}$ induces a $P_3\cup P_2$, a contradiction. So, $|N(x)\cap \{u_1,u_2,u_3,u_4\}|\ge3$. By symmetry, we may assume that $x\sim u_1, x\sim u_2$ and $x\sim u_3$. It is clear that all the vertices in $N_1$ which is complete to $\{u_1,u_2,u_3\}$ is independent as $G$ is $K_5$-free.

Suppose $N_1$ is not independent. We may by symmetry assume that $x'\in N_1$ such that $x'\sim x$. So, $x'$ is not complete to $\{u_1,u_2,u_3\}$, we may suppose that $u_1\not\in N(x')$, that means $N(x')\cap \{u_1,u_2,u_3,u_4\}=\{u_2,u_3,u_4\}$. Then $\{x,x',u_1,u_2,u_4\}$ induces a $W_4$, a contradiction. Therefore, $N_1$ is independent. With the same arguments, we have that $N_2$ and $N_4$ are both independent.

Suppose $N_1\cup N_2$ is not independent. We may by symmetry assume $y\in N_2$ such that $y\sim x$. Then $y$ is not complete to $\{u_1,u_2,u_3\}$, we may suppose that $u_1\not\in N(y)$, that means $N(y)\cap \{u_1,u_2,u_3,u_4\}=\{u_2,u_3,u_4\}$. Then $\{x,y,u_1,u_2,u_4\}$ induces a $W_4$, a contradiction. Therefore, $N_1\cup N_2$ is independent. With the same arguments, we have that $N_1\cup N_4$ and $N_2\cup N_4$ are both independent. This proves (\ref{eqa-3}).

Since $G$ is $(P_3\cup P_2, W_4)$-free, it follows that $G[N_3\setminus N_4]$ is $(P_3\cup P_2, C_4)$-free. Moreover, $\omega(G[N_3\setminus N_4])\le2$ as $G$ is $K_5$-free. By Lemma \ref{C4}, we have that $\chi(G[N_3\setminus N_4])\le 3$. Obviously, $G[M(\{v_1,v_2\})]$ is $P_3$-free, which implies that $\chi(G[M[\{v_1,v_2\}]])\le4$. Therefore, $\chi(G)\le \chi (G[N(\{v_1,v_2\})])+\chi(G[M[\{v_1,v_2\}]])\le \chi(G[N_1\cup N_2\cup N_4])+\chi(G[N_3\setminus N_4])+\chi(G[M[\{v_1,v_2\}]])\le1+3+4=8$. This proves Claim \ref{c4}.     \qed

\begin{claim}\label{c5}
	If $G$ is $(2K_3, P_2\cup K_4)$-free and contains an induced $P_2\cup K_3$, then $\chi(G)\le8$.
\end{claim}
\pf Let $G[\{v_1,v_2,u_1,u_2,u_3\}]$ be an induced $P_2\cup K_3$ in $G$ such that $u_1u_2u_3u_1$ is a $K_3$.
Let $v\in N(\{v_1,v_2\})$. Suppose that $N(v)\cap \{u_1,u_2,u_3\}=\emptyset$. Then $\{v_1,v_2,v,u_1,u_2\}$ induces a $P_3\cup P_2$ if $v\not\in N(v_1)\cap N(v_2)$, and $\{v_1,v_2,v,u_1,u_2,u_3\}$ induces a $2K_3$ if $v\in N(v_1)\cap N(v_2)$, both are contradictions. So, $N(v)\cap \{u_1,u_2,u_3\}\ne\emptyset$. Let $N_i=\{v\in N(\{v_1,v_2\}) | |N(v)\cap \{u_1,u_2,u_3\}|=i\}$, for $i\in\{1,2,3\}$. Note that $N(\{v_1,v_2\})=\cup_{i=1}^3 N_i$.

It is obvious that $N_1\subseteq N(v_1)\cap N(v_2)$ to forbid an induced $P_3\cup P_2$. Suppose $N_1$ has two vertices $x_1$ and $x_2$ such that $N(x_1)\cap \{u_1,u_2,u_3\}=N(x_2)\cap \{u_1,u_2,u_3\}=\{u_1\}$. Then $\{v_1,v_2,x_1,x_2,u_2,u_3\}$ induces a $P_2\cup K_4$ if $x\sim x_2$, and $\{x_1,v_1,x_2,u_2,u_3\}$ induces a $P_3\cup P_2$ if $x_1\not\sim x_2$, both are contradictions. So, we may assume that $N_1=\{x_1,x_2,x_3\}$ such that $x_i\cap\{u_1,u_2,u_3\}=u_i$, for $i\in\{1,2,3\}$. (It is possible that $x_1$ or $x_2$ or $x_3$ does not exist.)

Since $G$ is $(P_3\cup P_2, 2K_3)$-free, we have that each component  of $G[M(\{v_1,v_2\})\setminus\{u_1,u_2,u_3\}]$ is an edge or a vertex. If $x_1x_2x_3x_1$ is a $K_3$, then $\{v_1,v_2,x_1,x_2,x_3\}$ induces a $K_5$, this contradiction implies that $G[\{x_1,x_2,x_3,u_1\}]$ is a bipartite graph. It is obvious that $G[M[\{v_1,v_2\}]\setminus\{u_1\}]$ is a bipartite graph. Therefore, we have that $\chi(G[N_1\cup M[\{v_1,v_2\}]])\le4$.

Let $y_1y_2$ be an edge of $\chi(G[N_2])$. Suppose $N(y_1)\cap \{u_1,u_2,u_3\}\ne N(y_2)\cap \{u_1,u_2,u_3\}$ that $N(y_1)\cap \{u_1,u_2,u_3\}=\{u_1,u_2\}$ and $N(y_2)\cap \{u_1,u_2,u_3\}=\{u_2,u_3\}$. Then $\{y_1,y_2,u_1,u_2,u_3\}$ induces a $W_4$, a contradiction. So, $N(y_1)\cap \{u_1,u_2,u_3\}= N(y_2)\cap \{u_1,u_2,u_3\}$. This implies that $G[N_2]$ is $(K_3,C_4)$-free as $G$ is $(K_5, W_4)$-free. By Lemma \ref{C4}, we have that $\chi(G[N_2])\le3$.

Note that $N_3$ is independent as $G$ is $K_5$-free. It follows that $\chi(G)\le \chi(G[N_1\cup M[\{v_1,v_2\}]])+\chi(G[N_2])+\chi(G[N_3])\le4+3+1=8$. This proves Claim \ref{c5}.
\qed

\medskip

\begin{figure}[htbp]\label{fig-3}
	\begin{center}
		\includegraphics[width=5cm]{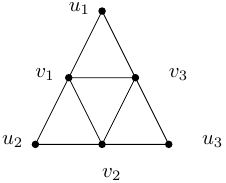}
	\end{center}
	\vskip -15pt
	\caption{Illustration of 4-triangle.}
\end{figure}

Let  $C=u_1v_2u_2v_2u_3v_3u_1$ be a 6-hole. We use {\em 4-triangle} (see Figure \ref{fig-3}) to denote a graph obtained from $C$ by connecting edges $v_1v_2$, $v_2v_3$ and $v_3v_1$.  

\begin{claim}\label{c6}
	If $G$ is $(P_2\cup K_3)$-free and contains an induced 4-triangle, then $\chi(G)\le8$.
\end{claim}
\pf Let $Q$ be an induced 4-triangle in $G$ with the same labels as in Figure 3, and $S=\{v_1,v_2,v_3\}$. Let $N_i=\{v\in V(G)\setminus S | |d(v)\cap S|=i\}$, for $i\in\{0,1,2,3\}$. Note that $N(S)=\cup_{i=1}^3 N_i$, and $M(S)=N_0$.

Let $A_i=\{v\in V(G)\setminus S | v $ is anticomplete to $S\setminus\{v_i\}\}$. We see that $N_0\cup N_1=\cup_{i=1}^3 A_i$.
Since $G$ is $(P_3\cup P_2, P_2\cup K_3)$-free, it follows that each component of $G[A_i]$ is a vertex or an edge for each $i\in\{1,2,3\}$.  For each $i\in\{1,2,3\}$, by the definition of $A_i$, we have that $A_i\cup S\setminus\{v_i\}$ induces a bipartite graph. So, $\chi(G[N_0\cup N_1\cup S])\le6$.

For $i\in\{1,2,3\}$, let $B_i=\{v\in N_2 | v \not\sim v_i\}$, we see that $N_2=\cup_{i=1}^3 B_i$, and $B_i\ne\emptyset$ as the $Q$ exists. Suppose $N_2$ has two adjacent vertices, say $x_1$ and $x_2$. Without loss of generality, we suppose that $x_1=u_1$. It is obvious that $x_2\not\in\{u_2,u_3\}$.  

Suppose $x_2\in B_2$, that is to say, $x_2\sim v_1, x_2\sim v_3$ and $x_2\not\sim v_2$. Then $x_2\not\sim u_2$ and $x_2\not\sim u_3$ as otherwise $\{x_2,v_1,u_2,v_2,v_3\}$ or $\{x_2,v_1,v_2,u_3,v_3\}$ induces a $W_4$, a contradiction. Now, $\{u_1,x_2,u_2,v_2,u_3\}$ induces a $P_3\cup P_2$, a contradiction. Therefore, $x_2\not\in B_2$. By symmetry, we may assume that $x_2\in B_1$, that is to say, $x_2\sim v_2, x_2\sim v_3$ and $x_2\not\sim v_1$.  But now, $\{u_1,v_1,v_2,x_2,v_3\}$ induces a $W_4$, a contradiction. So, $N_2$ is independent.

Note that $N_3$ is independent as $G$ is $K_5$-free. It follows that $\chi(G)\le \chi(G[N_0\cup N_1\cup S ])+\chi(G[N_2])+\chi(G[N_3])\le6+1+1=8$. This proves Claim \ref{c6}.
\qed

\begin{claim}\label{c7}
	If $G$ is $($$P_2\cup K_3$, 4-triangle$)$-free and contains an induced gem, then $\chi(G)\le8$.
\end{claim}

\pf Let $Q$ be an induced gem in $G$ with $V(Q)=\{v,u_1,u_2,u_3,u_4\}$ and $d_Q(v)=4$. Let $x\in M(v)$. Suppose $x\sim u_2$ and $x\sim u_3$. Then $\{v,u_1,u_2,u_3,x\}$ induces a $W_4$ if $x\sim u_1$, and $\{v,u_2,u_3,u_4,x\}$ induces a $W_4$ if $x\sim u_4$, and $\{v,u_1,u_2,u_3,u_4,x\}$ induces a 4-triangle if $x$ is anticomplete to $\{u_1,u_4\}$, all are contradictions. So, $x\not\sim u_2$ or $x\not\sim u_3$.

Let $M=\{x\in M(v) | x\not\sim u_2\}$ and $M'=\{x\in M(v) | x\not\sim u_3\}$. Then $G[M]$ is $(P_3, K_3)$-free as $M$ is anticomplete to $\{v,u_2\}$ and $G$ is $(P_3\cup P_2, P_2\cup  K_3)$-free. So, $G[M]$ is bipartite. By symmetry, $G[M']$ is bipartite.

Since $G$ is $(P_3\cup P_2, W_4, K_5)$-free, we have that $G[N(v)]$ is $(P_3\cup P_2, C_4, K_4)$-free, which implies that $\chi(G[N(v)])\le4$ by Lemma \ref{C4}.  So $\chi(G)\le \chi(G[N(v)])+\chi(G[M[v]])\le \chi(G[N(v)])+\chi(G[M])+\chi(G[M'])\le4+2+2=8$. This proves Claim \ref{c7}.\qed

\medskip

By Claim \ref{c3}, \ref{c4}, \ref{c5}, \ref{c6}, \ref{c7}, we may assume that $G$ is ($P_3\cup P_2$, gem)-free. By Theorem \ref{gem}, we have that $\chi(G)\le8$. This proves Lemma \ref{=4}.\qed

\section{$K_4$-free}

In this section, we consider ($P_3\cup P_2, W_4, K_4$)-free graphs. Let $C=v_1v_2u_2u_1u_3v_3v_1$ be a 6-hole. We use $F_1$ to denote a graph obtained from $C$ by connecting edges $v_2v_3$ and $u_2u_3$,  use $F_3$ to denote a graph obtained from $C$ by connecting edges $v_2v_3$, use $F_4$ to denote a graph obtained from $C$ by connecting edges $v_2v_3,u_2u_3$ and $u_2v_3$ (see Figure 4).

Let $C'=v_1v_2u_3u_1u_2v_1$ be a 5-hole. We use $F_2$ to denote a graph obtained from $C$ by adding a vertex $u_4$ such that $u_4$ is complete to $\{u_1,u_2,u_3\}$ and anticomplete to $\{v_1,v_2\}$ (see Figure 4).

A {\em hammer} is the graph obtained by identifying one vertex of a $K_3$ and one end vertex of a $P_3$ (see Figure 4).

\begin{figure}[htbp]\label{fig-4}
	\begin{center}
		\includegraphics[width=17cm]{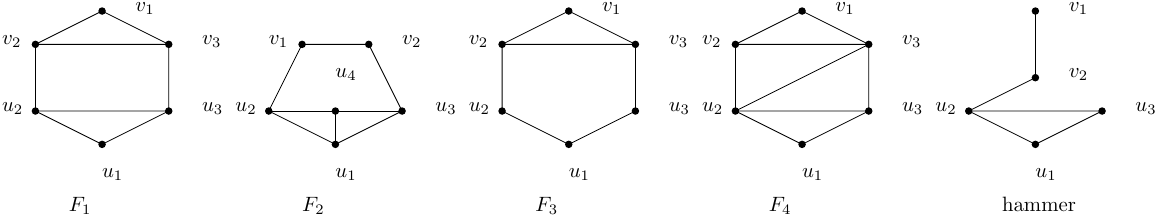}
	\end{center}
	\vskip -15pt
	\caption{Illustration of $F_1,F_2,F_3,F_4$ and hammer.}
\end{figure}

\begin{lemma}\label{=3}
	If $G$ is a $(P_3\cup P_2, W_4, K_4)$-free graph, then $\chi(G)\le6$.
\end{lemma}

\pf Let $G$ be a $(P_3\cup P_2, W_4, K_4)$-free graph, we will complete the proof of Lemma \ref{=3} by the following six claims. 

\begin{claim}\label{c4-1}
	If $G$ contains an induced $P_2\cup K_3$, then $\chi(G)\le6$.
\end{claim}
\pf Let $G[\{v_1,v_2,u_1,u_2,u_3\}]$ be an induced $P_2\cup K_3$ in $G$ such that $u_1u_2u_3u_1$ is a $K_3$.
Let $N= N(\{v_1,v_2\})$. Then no vertex of $N$ is complete to $\{u_1,u_2,u_3\}$ as $G$ is $K_4$-free. Let $N_1=\{v | N(v)\cap \{v_1,v_2\}=\{v_1\}\}, N_2=\{v | N(v)\cap \{v_1,v_2\}=\{v_2\}\}$, and $N_3=\{v | \{v_1,v_2\}\subseteq  N(v)\}$. Notice that $N=\cup_{i=1}^3 N_i$. 

Let $v$ be any vertex in $N_1$. If $|\{u_1,u_2,u_3\}\setminus N(v)|\ge2$, let $v\not\sim u_1$ and $v\not\sim u_2$, then $\{v,v_1,v_2,u_1,u_2\}$ induces a $P_3\cup P_2$, a contradiction. So, $|N(v)\cap \{u_1,u_2,u_3\}|=2$. Suppose there exists two adjacent vertices $u$ and $u'$ in $N_1$. Then $|N(u)\cap \{u_1,u_2,u_3\}|=|N(u')\cap \{u_1,u_2,u_3\}|=2$. We may by symmetry assume that $N(u)\cap \{u_1,u_2,u_3\}=\{u_1,u_2\}$. Then $N(u')\cap \{u_1,u_2,u_3\}\ne\{u_1,u_2\}$ as $G$ is $K_4$-free. So, we may assume that $N(u')\cap \{u_1,u_2,u_3\}=\{u_1,u_3\}$. But now, $\{u,u',u_1,u_2,u_3\}$ induces a $W_4$, a contradiction. Therefore, $N_1$ is independent, and by symmetry, $N_2$ is independent.

It is obvious that $N_3$ is independent as $G$ is $K_4$-free. So, $\chi(G)\le\chi(G[N])+\chi(G[M[\{v_1,v_2\}]])\le3+3=6$. This proves Claim \ref{c4-1}. \qed

\begin{claim}\label{c4-2}
	If $G$ is $(P_2\cup K_3)$-free and contains an induced $F_1$, then $\chi(G)\le6$.
\end{claim}
\pf Let $Q$ be an induced $F_1$ in $G$ with the same labels as in Figure 4, and $S=\{v_1,v_2,v_3\}$.  Let $N_i=\{v\in V(G)\setminus S | |d(v)\cap S|=i\}$, for $i\in\{0,1,2\}$. Note that $N(S)=\cup_{i=1}^2N_i$ as $G$ is $K_4$-free, and $M(S)=N_0$.

Let $A_i=\{v\in  N(S) | v $ is anticomplete to $S\setminus\{v_i\}\}$. We see that $N_1=\cup_{i=1}^3 A_i$.
Let $v$ be any vertex in $A_1$. Then either $v\sim u_1$ or $v$ is complete to $\{u_2,u_3\}$. Otherwise, we may by symmetry assume that $v\not\sim u_1$ and $v\not\sim u_2$, which implies that $\{v,v_1,v_3,u_1,u_2\}$ induces a $P_3\cup P_2$, a contradiction. Suppose there exist two adjacent vertices $u$ and $u'$ in $A_1$. If $u$ and $u'$ are both adjacent to $u_1$, then $\{u,u',u_1,v_2,v_3\}$ induces a $P_2\cup K_3$, a contradiction.  So, we may by symmetry assume that $u\not\sim u_1, u\sim u_2$ and $u\sim u_3$. But now, $\{u,u',u_1,v_2,v_3\}$ induces a $P_3\cup P_2$ if $u'\sim u_1$ and $\{u,u',u_2,u_3\}$ induces a $K_4$ if $u'\not\sim u_1$, both are contradictions. Therefore, $A_1$ is independent, and by symmetry, $A_2$ and $A_3$ are both independent.

Suppose there exist two adjacent vertices $x$ and $x'$ in $N_2$. Then $|N(x)\cap \{v_1,v_2,v_3\}|=|N(x')\cap \{v_1,v_2,v_3\}|=2$. We may by symmetry assume that $N(x)\cap \{v_1,v_2,v_3\}=\{v_1,v_2\}$. Then $N(x')\cap \{v_1,v_2,v_3\}\ne\{v_1,v_2\}$ as $G$ is $K_4$-free. So, we may assume that $N(x')\cap \{v_1,v_2,v_3\}=\{v_1,v_3\}$. But now, $\{x,x',v_1,v_2,v_3\}$ induces a $W_4$, a contradiction. Therefore, $N_2$ is independent.

Since $v_i$ is anticomplete to $A_i$ for each $i\in\{1,2,3\}$, we have that $\chi(G[N_1\cup\{v_1,v_2,v_3\}])\le3$. Moreover, $\chi(G[N_0])\le1$ as $G$ is $P_2\cup K_3$-free. So, $\chi(G)\le \chi(G[ N_1\cup\{v_1,v_2,v_3\}])+\chi(G[N_2])+\chi(G[N_0])\le3+1+1=5<6$. This proves Claim \ref{c4-2}. \qed

\begin{claim}\label{c4-3}
	If $G$ is $(P_2\cup K_3, F_1)$-free and contains an induced $F_2$, then $\chi(G)\le6$.
\end{claim}

\pf We will firstly prove that 

\begin{equation}\label{eqa4-1}
	\mbox{for each vertex $v\in V(G)$, $G[N(v)]$ is bipartite if $G[M(v)]$ contains a $K_3$}.
\end{equation}

Let $v\in V(G)$, and $K$ be a $K_3$ in $G[M(v)]$ with $V(K)=\{t_1,t_2,t_3\}$. Let $N_i=\{u\in N(v)| |d_K(u)|=i\}$, for $i\in\{1,2\}$. Note that $N(v)=\cup_{i=1}^2N_i$ as $G$ is $(K_4, P_2\cup K_3)$-free.

Suppose there exist two adjacent vertices $x$ and $x'$ in $N_1$. Then $|N_K(x)|=|N_K(x')|=1$. If $N_K(x)=N_K(x')$, say $t_1$, then $\{x,x',v,t_2,t_3\}$ induces a $P_2\cup K_3$, a contradiction. So, we may suppose that $N_K(x)\ne N_K(x')$. We may assume that $x\sim t_1$ and $x'\sim t_2$. Then $\{x,x',v,t_1,t_2,t_3\}$ induces a $F_1$, a contradiction. So, $N_1$ is independent.

Suppose there exist two adjacent vertices $y$ and $y'$ in $N_2$. Then $|N_K(y)|=|N_K(y')|=2$. If $N_K(y)=N_K(y')$, say $\{t_1,t_2\}$, then $\{x,x',t_1,t_2\}$ induces a $K_4$, a contradiction. So, $N_K(y)\ne N_K(y')$, which implies that $\{y,y',t_1,t_2,t_3\}$ induces a $W_4$, a contradiction. 

So, $N_1$ and $N_2$ are both independent. This proves (\ref{eqa4-1}).

Let $Q$ be an induced $F_2$ in $G$ with the same labels as in Figure 4. Since both $G[M(v_1)]$ and $G[M(v_2)]$ contain a $K_3$, we have that both $G[N(v_1)]$ and $G[N(v_2)]$ are bipartite. Consequently, $G[M(\{v_1,v_2\})]$ is bipartite as $G$ is $(P_3\cup P_2, P_2\cup K_3)$-free. So, $\chi(G)\le \chi(G[N(\{v_1,v_2\})])+\chi(G[M[\{v_1,v_2\}]])\le \chi(G[N(v_1)])+\chi(G[N(v_2)])+\chi(G[M[\{v_1,v_2\}]])\le2+2+2=6$. This proves Claim \ref{c4-3}. \qed

\begin{claim}\label{c4-4}
	If $G$ is $(P_2\cup K_3, F_1, F_2)$-free and contains an induced $F_3$, then $\chi(G)\le6$.
\end{claim}
\pf Let $Q$ be an induced $F_3$ in $G$ with the same labels as in Figure 4, and $S=\{v_1,v_2,v_3\}$.  Let $N_i=\{v\in V(G)\setminus S | |d(v)\cap S|=i\}$, for $i\in\{0,1,2\}$. Note that $N(S)=\cup_{i=1}^2N_i$ as $G$ is $K_4$-free, and $M(S)=N_0$.

Let $A_i=\{v\in N(S) | v $ is anticomplete to $S\setminus\{v_i\}\}$. We see that $N_1=\cup_{i=1}^3 A_i$.
By the same arguments in the proof of Claim \ref{c4-2}, we only need to prove that $A_1$ is independent.

Let $v$ be any vertex in $A_1$. Then either $v\sim u_1$ or $v$ is complete to $\{u_2,u_3\}$. Otherwise, we may by symmetry assume that $v\not\sim u_1$ and $v\not\sim u_2$, which implies that $\{v,v_1,v_3,u_1,u_2\}$ induces a $P_3\cup P_2$, a contradiction. Suppose there exist two adjacent vertices $u$ and $u'$ in $A_1$. If $u$ and $u'$ are both adjacent to $u_1$, then $\{u,u',u_1,v_2,v_3\}$ induces a $P_2\cup K_3$, a contradiction.  So, we may by symmetry assume that $u\not\sim u_1, u\sim u_2$ and $u\sim u_3$. But now, $\{u,u',u_1,v_2,v_3\}$ induces a $P_3\cup P_2$ if $u'\sim u_1$ and $\{u,u',v_2,v_3,u_2,u_3\}$ induces a $F_2$ if $u'\not\sim u_1$, both are contradictions. Therefore, $A_1$ is independent.

This proves Claim \ref{c4-4}. \qed

\begin{claim}\label{c4-5}
	If $G$ is $(P_2\cup K_3, F_1, F_2, F_3)$-free and contains an induced $F_4$, then $\chi(G)\le6$.
\end{claim}
\pf Let $Q$ be an induced $F_4$ in $G$ with the same labels as in Figure 4, and $S=\{v_1,v_2,v_3\}$.  Let $N_i=\{v\in V(G)\setminus S | |d(v)\cap S|=i\}$, for $i\in\{0,1,2\}$. Note that $N(S)=\cup_{i=1}^2N_i$ as $G$ is $K_4$-free, and $M(S)=N_0$.

Let $A_i=\{v\in N(S) | v $ is anticomplete to $S\setminus\{v_i\}\}$. We see that $N_1=\cup_{i=1}^3 A_i$.
By the same arguments in the proof of Claim \ref{c4-2}, we only need to prove that $A_1$ is independent.

Let $v$ be any vertex in $A_1$. Then either $v\sim u_1$ or $v$ is complete to $\{u_2,u_3\}$. Otherwise, we may assume that $v\not\sim u_1$. If $v\not\sim u_3$, then $\{v,v_1,v_2,u_1,u_3\}$ induces a $P_3\cup P_2$, a contradiction. So, $v\sim u_3$ and $v\not\sim v_2$. Then $\{v,v_1,v_2,u_1,u_2,u_3\}$ induces a $F_3$, a contradictions. Therefore, either $v\sim u_1$ or $v$ is complete to $\{u_2,u_3\}$.

Suppose there exist two adjacent vertices $u$ and $u'$ in $A_1$. If $u$ and $u'$ are both adjacent to $u_1$, then $\{u,u',u_1,v_2,v_3\}$ induces a $P_2\cup K_3$, a contradiction.  So, we may by symmetry assume that $u\not\sim u_1, u\sim u_2$ and $u\sim u_3$. But now, $\{u,u',u_1,v_2,v_3\}$ induces a $P_3\cup P_2$ if $u'\sim u_1$ and $\{u,u',u_2,u_3\}$ induces a $K_4$ if $u'\not\sim u_1$, both are contradictions. Therefore, $A_1$ is independent.

This proves Claim \ref{c4-5}. \qed

\begin{claim}\label{c4-6}
	If $G$ is $(P_2\cup K_3, F_1, F_2, F_3, F_4)$-free and contains an induced hammer, then $\chi(G)\le6$.
\end{claim}
\pf Let $Q$ be an induced hammer in $G$ with the same labels as in Figure 4. Since $G[N(v_2)]$ is $C_4$-free and $\omega(G[N(v_2)])\le2$, we have that $\chi(G[N(v_2)])\le3$ by Theorem \ref{C4}. 

Let $N_1=N(\{v_1,v_2\})\setminus N(v_2)$. Suppose there exist two adjacent vertices $u$ and $u'$ in $N_1$. Note that neither $u$ nor $u'$ is complete to $\{u_1,u_2,u_3\}$ as $G$ is $K_4$-free. If $u$ is anticomplete to $\{u_1,u_3\}$, then $\{v_2,v_1,u,u_1,u_3\}$ induces a $P_2\cup K_3$, a contradiction. So, $N(u)\cap\{u_1,u_3\}\ne\emptyset$. By symmetry, $N(u')\cap\{u_1,u_3\}\ne\emptyset$. 

Suppose $u\sim u_3$. If $u$ is anticomplete to $\{u_1,u_2\}$, then $\{u_1,u_2,u_3,u,v_1,v_2\}$ induces a $F_3$, a contradiction. So, we firstly suppose $u\sim u_2$. If $u'\sim u_3$, then $\{v_1,u,u',u_1,u_2,u_3\}$ induces a $F_4$ when $u'$ is anticomplete to $\{u_1,u_2\}$, and $\{u',u_1,u_2,u_3,v_1,v_2\}$ induces a $F_2$ when $u'\sim u_1$ and $u'\not\sim u_2$, and $\{u,u',u_2,u_3\}$ induces a $K_4$ when $u'\sim u_2$, all are contradictions. So, $u'\not\sim u_3$ and $u'\sim u_1$. Then $\{u_1,u_2,u_3,u',v_1,v_2\}$ induces a $F_3$ if $u'\not\sim u_2$, and $\{u_2,u,u',u_1,u_3\}$ induces a $W_4$ if $u'\sim u_2$, both are contradictions. Thus $u\not\sim u_2$ and $u\sim u_1$. Then $u'$ is not complete to $\{u_1,u_3\}$ to forbid a $K_4$ on $\{u,u',u_1,u_3\}$. If $u'\sim u_3$, then $\{v_1,u,u',u_1,u_2,u_3\}$ induces a $F_4$ when $u'\not\sim u_2$, and $\{u,u',u_1,u_2,u_3\}$ induces a $W_4$ when $u'\sim u_2$, both are contradictions. Therefore, $u\not\sim u_3$.

So, $u\sim u_1$. Then $u\sim u_2$ as otherwise $\{v_1,u',u_1,u_2,u_3\}$ induces a $P_2\cup K_3$, a contradiction. If $u'$ is complete to $\{u_1,u_3\}$, then $\{u_1,u_2,u_3,u,u'\}$ induces a $W_4$, a contradiction. 
If $u'\sim u_1$ and $u'\not\sim u_3$ , then $\{u_1,u_2,u_3,u,u',v_1\}$ induces a $F_4$, a contradiction. If $u'\sim u_3$ and $u'\not\sim u_1$ , then $\{u_1,u_2,u_3,u',v_1,v_2\}$ induces a $F_3$, a contradiction.  Therefore, $u\not\sim u_1$. This proves that $N_1$ is independent.

Since $M(\{v_1,v_2\})$ is $(P_3,K_3)$-free, we have that each component of $G[M(\{v_1,v_2\})]$ is a vertex or an edge. So, $\chi(G)\le \chi(G[N(\{v_1,v_2\})])+\chi(G[M[\{v_1,v_2\}]])\le \chi(G[N_1])+\chi(G[N(v_2)])+\chi(G[M[\{v_1,v_2\}]])\le1+3+2=6$. This proves Claim \ref{c4-6}. \qed

\medskip

By Claim \ref{c4-1}, \ref{c4-2}, \ref{c4-3}, \ref{c4-4}, \ref{c4-5} and \ref{c4-6}, we may assume that $G$ is $($$P_2\cup K_3, F_1, F_2, F_3, F_4$, hammer$)$-free. 

Suppose $G$ contains an induced $K_1\cup K_3$. Let $Q$ be an induced $K_1\cup K_3$ in $G$ such that $d_Q(v)=0$. Since $G$ is $($$P_2\cup K_3, W_4,K_4$, hammer$)$-free, we can easily verify that $G[N(v)]$ is independent. Let $v'\in N(v)$. By Theorem \ref{C4}, $\chi(G[N(v')]) \le3$ as $G[N(v')]$ is $C_4$-free and $\omega(G[N(v')])\le2$. Moreover, since $G[M(\{v,v'\})]$ is $(P_3,K_3)$-free, we have that each component of $G[M(\{v,v'\})]$ is a vertex or an edge. So, $\chi(G)\le \chi(G[N(\{v,v'\})])+\chi(G[M[\{v,v'\}]])\le \chi(G[N(v)])+\chi(G[N(v')])+\chi(G[M[\{v_1,v_2\}]])\le1+3+2=6$. 

So, we may assume that $G$ is $(K_1\cup K_3)$-free. Consequently, we may suppose that $G$ contains a $K_3$ by Theorem \ref{C3}. Now, by Theorem \ref{K3}, we have that $\chi(G)\le6$. This proves Lemma \ref{=3}.\qed

\end{document}